\documentclass{article}

\title{A generalization of Kummer's identity}

\author{Raimundas Vid\=unas\footnote{Supported by NWO,
project number 613-06-565.}\\ \em University of Amsterdam}

\newtheorem{theorem}{Theorem}

\newcommand{\hpg}[5]{{}_{#1}\mbox{\rm F}_{\!#2}\!
  \left(\left.{#3 \atop #4}\right| #5 \right) }
\newcommand{\hpgg}[4]{{}_{#1}\mbox{\rm F}_{\!#2}\!
  \left({#3 \atop #4}\right) }
\newcommand{\hpgo}[2]{{}_{#1}\mbox{\rm F}_{\!#2}}

\newcommand{\poh}[2]{(#1)_{#2}}

\newcommand{\re}[1]{\mbox{Re}(#1)}

\newcommand{\CC}{\mbox{\bf C}}

\begin{document}

\maketitle

\begin{abstract}
The well-known Kummer's formula evaluates 
the hypergeometric series $\hpg{2}{1}{A,\,B}{C}{-1}$ 
when relation $C\!-\!A\!+\!B\!=\!1$ holds. This paper deals with evaluation
of $\hpgo{2}{1}(-1)$ series in the case when $C\!-\!A\!+\!B$ is an integer.
Such a series is expressed as a sum of two $\Gamma$-terms multiplied by
terminating $\hpgo{3}{2}(1)$ series. A few such formulas were essentially
known to Whipple in 1920's. Here we give a simpler and more complete
overview of this type of evaluations. Additionally, algorithmic aspects
of evaluating hypergeometric series are considered. We illustrate
Zeilberger's method and discuss its applicability to non-terminating series,
and present a couple of similar generalizations of other known formulas.

\end{abstract}

\section{The generalization}

The subject of this paper is a generalization of Kummer's identity
(see \cite{kummer}, \cite[2.3]{bailey}, or Cor.~3.1.2 in \cite{specfaar}):
\begin{equation} \label{ikummer}
\hpg{2}{1}{a,\;b}{1+a-b}{\,-1}=\frac{\Gamma(1+a-b)\;\Gamma(1+\frac{a}{2})}
{\Gamma(1+a)\;\Gamma(1+\frac{a}{2}-b)}. 
\end{equation}
The hypergeometric series on the left is defined if $a-b$ is not a negative
integer, and it is absolutely convergent for $\re{b}<1/2$. 
After analytic continuation of $\hpg{2}{1}{a,\,b}{\!1+a-b}{z}$ on
$\CC\setminus [1,\infty)$ and after division of both sides by
$\Gamma(1\!+\!a\!-\!b)$ the formula has meaning and is correct for
all complex $a,b$. 
In this paper, whenever $\hpg{2}{1}{\!A,\,B}{C}{z}$ denotes a well-defined
hypergeometric series, it also denotes its analytic continuation on
$\CC\setminus [1,\infty)$.

The generalization to be considered evaluates the hypergeometric series
$\hpg{2}{1}{\!A,\,B}{C}{-1}$ whenever $C\!-\!A\!+\!B$ is any integer.
In the terminology of \cite{specfaar}, our generalization applies
to $\hpgo{2}{1}(-1)$ series which are {\em contiguous} to a series for
Kummer's formula (\ref{ikummer}).
As it is known (see \cite[2.5]{specfaar}), the 15 classical Gauss contiguity
relations can be iterated to produce a linear relation between any three
contiguous $\hpgo{2}{1}(z)$ series, with coefficients being rational functions
in the parameters of those series. This also applies to their analytic
extensions. The generalized formula is such a relation in explicit form
between contiguous $\hpg{2}{1}{\!a+n,\,b}{a-b}{-\!1}$,
$\hpg{2}{1}{a,\,b}{\!1+a-b}{-\!1}$ and $\hpg{2}{1}{\!a-1,\,b}{a-b}{-\!1}$,
where $n$ is an integer, and the last two series are evaluated using Kummer's
identity (\ref{ikummer}). The coefficient to the first series cannot be the
zero function because the quotient of the other two series is not in
$\CC(a,b,n)$. In the generalized formula these coefficients are written
as terminating $\hpgo{3}{2}(1)$ series.

We write the generalization in the form
\begin{equation} \label{genkum}
\hpg{2}{1}{a+n,\,b}{a-b}{-1}= P(n)\,
\,\frac{\Gamma(a-b)\,\Gamma(\frac{a+1}{2})}
{\Gamma(a)\,\Gamma(\frac{a+1}{2}-b)}+Q(n)\,
\,\frac{\Gamma(a-b)\,\Gamma(\frac{a}{2})}
{\Gamma(a)\,\Gamma(\frac{a}{2}-b)}.
\end{equation}
Here the two $\Gamma$-terms are equal to $\hpg{2}{1}{\!a-1,\,b}{a-b}{-\!1}$
and $\frac{a-b}{a-2b}\,\hpg{2}{1}{\!a,\,b}{\!1+a-b}{-\!1}$ respectively,
and $P(n)$, $Q(n)$ are rational functions in $a,b$ for every integer $n$.
The most convenient expressions for $P(n)$ and $Q(n)$ are summarized in
the three theorems below. 
In fact, expressions of $\hpgo{2}{1}(-1)$ series in (\ref{genkum}) in
terms of terminating series and $\Gamma$-function were known to Whipple,
see \cite{whipple30}. His formulas (8.3) and (8.41) would express the
$\hpgo{2}{1}(-1)$ series in (\ref{genkum}) in terminating series for
negative or positive $n$, respectively. Whipple's formulas (11.5,51) form
the statement of Theorem \ref{gtheorem1} below. Whipple derived them as
a consequence of  transformations of $\hpgo{3}{2}(1)$ series allied to
general $\hpgo{2}{1}(-1)$ series, and from \cite[(2.6,7)]{fox27} where
some $\hpgo{2}{1}(1/2)$ series are expressed in terms of terminating series.
However, Whipple's main concern was the relations of general $\hpgo{2}{1}(-1)$
and $\hpgo{3}{2}(1)$ series. As we will see, his approach is not convenient
when some of those series terminate.

In this paper we strive for a clear overview of possible expressions
for $P(n)$ and $Q(n)$ in terms of terminating $\hpgo{3}{2}(1)$ series,
with simpler proofs. 
Another aim is to consider algorithmic aspects of evaluating hypergeometric
series. In particular, we specialize formula (\ref{genkum}) to
two-term identities, which  however seem to be beyond Zeilberger's approach.
Also a few evaluations similar to (\ref{genkum}) are presented.
Specifically, we evaluate hypergeometric series which are contiguous
to the $\hpgo{2}{1}(1/4)$ and $\hpgo{3}{2}(1)$ series in Gosper's and
Dixon's identities, see (\ref{gengosper}-\ref{gendixon}).

In the following theorems we summarize the most convenient expressions for
$P(n)$ and $Q(n)$. A few more such expressions are presented
in (\ref{altpn}-\ref{altqn1}).
\begin{theorem} \label{gtheorem1}
Suppose that $n$ is a non-negative integer {\rm(}or $-1${\rm)},
and $a,b$ are complex numbers such that $(a)_n\neq 0$ and $a-b$ is not zero
or a negative integer. Then the coefficients $P(n)$ and $Q(n)$ in formula
{\rm(\ref{genkum})}
can be written as:
\begin{eqnarray} \label{whippn}
P(n) & = & \frac{1}{2^{n+1}}\;
\hpgg{3}{2}{-\frac{n}{2},\,-\frac{n+1}{2},\,\frac{a}{2}-b}
{\frac{1}{2},\;\frac{a}{2}}, \\
\label{whipqn} Q(n) & = & \frac{n\!+\!1}{2^{n+1}}\;
\hpgg{3}{2}{-\frac{n-1}{2},\,-\frac{n}{2},\,\frac{a+1}{2}-b}
{\frac{3}{2},\;\frac{a+1}{2}}.
\end{eqnarray}
\end{theorem}
\begin{theorem} \label{gtheorem2}
Suppose that $n$ is a non-negative integer, and $a,b$ are complex such
that $(a)_n\neq 0$, and $a-b$ is not zero or a negative integer. Then the
coefficients $P(n)$ and $Q(n)$ in formula {\rm(\ref{genkum})} can be
written as:
\begin{eqnarray} \label{pqdef1}
P(n) = \frac{1}{2}\, 
\hpgg{3}{2}{-\frac{n}{2},\,-\frac{n+1}{2},\,b}{-n,\;\frac{a}{2}}, 
\qquad  Q(n) = \frac{1}{2}\, 
\hpgg{3}{2}{-\frac{n-1}{2},\,-\frac{n}{2},\,b}{-n,\;\frac{a+1}{2}}.
\end{eqnarray}
The hypergeometric sums should be interpreted as terminating series with
{\rm(}up to $\pm 1${\rm)} $\lfloor n/2 \rfloor$ terms.
\end{theorem}
\begin{theorem} \label{gtheoremneg}
Let $P(n,a,b)$ and $Q(n,a,b)$ denote the coefficients $P(n)$ and $Q(n)$
in {\rm(\ref{genkum})} as functions of $a,b$ as well. If $n$ is a
non-negative integer, and $a,b\not\in\{0,1,\ldots,n\}$ then
\begin{eqnarray}
P(-n\!-\!1,\,a,\,b) & = & 2^{2n}\,\frac{(1\!-\!\frac{a}{2})_n}{(1\!-\!b)_n}
\;P(n\!-\!1,\,a\!-\!2n,\,b\!-\!n),\\
Q(-n\!-\!1,\,a,\,b) & = & \!-2^{2n}\,\frac{(\frac{1-a}{2})_n}{(1\!-\!b)_n}
\;Q(n\!-\!1,\,a\!-\!2n,\,b\!-\!n).
\end{eqnarray}
\end{theorem}
Because of the last theorem we do not give expressions for $P(n)$ and $Q(n)$
for a negative $n$, except (\ref{whippn1}-\ref{whipqn1}) in the proof of
Theorem \ref{gtheoremneg}.

These theorems are proved in section \ref{clproof}. There we also
overview transformations between other expressions for $P(n)$ and $Q(n)$,
and give a survey of Whipple's approach in \cite{whipple30}.
In section \ref{ctproof} Theorem \ref{gtheorem2} is proved using the more
universal Zeilberger's method. The key observation is that the sequences
$P(n)$ and $Q(n)$ satisfy the same recurrence relation as the left-hand side
of (\ref{genkum}). Theorem \ref{gtheorem1} can also be proved in this way.
Notice that any different expressions for $P(n)$ and $Q(n)$ must represent
the same rational functions in $a$, $b$ for every $n$, because the quotient
of the $\Gamma$-terms in (\ref{genkum}) is not in $\CC(a,b)$.
Section \ref{comalge} is devoted to algorithmic aspects of evaluation of
hypergeometric series, with similar generalizations of Dixon's and Gosper's
identities. 

{\em Acknowledgments.} The author would like to thank Richard Askey and
Tom Koornwinder for useful suggestions, in particular for ideas of the easy
proof(s) in section \ref{clproof}, and Dennis Stanton for the references to
Whipple.

\section{Classical proof}
\label{clproof}

We assume here $\re{a/2}>\re{b}>0$. One can simply check that Theorems
\ref{gtheorem1} and \ref{gtheorem2} must hold for the analytic continuation
of the $\hpgo{2}{1}(-1)$ series as well.

To prove Theorem \ref{gtheorem1} we recall
Whipple's identity \cite[(8.41)]{whipple30}
\begin{equation} \label{whip841}
\hpg{2}{1}{A,\,B}{C}{-1}=\frac{\Gamma(C)}{2\cdot\Gamma(A)}\,
\sum_{k=0}^{\infty} (-1)^k\,\frac{(C\!-\!A\!+\!B\!-\!1)_k}{k!}\,
\frac{\Gamma(\frac{A}{2}\!+\!\frac{k}{2})}
{\Gamma(C\!-\!\frac{A}{2}\!+\!\frac{k}{2})}.
\end{equation}
As it was communicated by Askey, this identity can be easily proved using
Euler's integral representation (\cite[2.12(1)]{bateman}) for the
$\hpgo{2}{1}(z)$ series. One has to rearrange the integrand as 
\begin{equation} \label{rearrange}
t^{A-1}\,(1\!-\!t)^{C-A-1}\,(1\!+\!t)^{-B}=
t^{A-1}\,(1\!+\!t)^{-C+A-B+1}\,(1\!-\!t^2)^{C-A-1},
\end{equation}
expand $(1\!+\!t)^{-C+A-B+1}$ as series, interchange integration and
summation, change the variable $t\!\mapsto\!\sqrt{s}$, and
recognize the beta-integral \cite[1.5(1)]{bateman}.

We apply\footnote{The same could be done directly to
$\hpg{2}{1}{a+n,\,b}{a-b}{-1}$, of course. We would get
less-convenient formula
\begin{eqnarray*}
\hpg{2}{1}{a+n,\,b}{a-b}{-1} & = & \frac{1}{2}\,
\frac{\Gamma(a\!-\!b)\,\Gamma(\frac{a+n}{2})}
{\Gamma(a\!+\!n)\,\Gamma(\frac{a-n}{2}\!-\!b)}\,
\hpgg{3}{2}{-\frac{n}{2},\,-\frac{n+1}{2},\,\frac{a+n}{2}}
{\frac{1}{2},\;\frac{a-n}{2}-b} \\ & & \hspace{-6pt}
+\frac{n\!+\!1}{2} \frac{\Gamma(a\!-\!b)\,\Gamma(\frac{a+n+1}{2})}
{\Gamma(a\!+\!n) \Gamma(\frac{a-n+1}{2}\!-\!b)}
\hpgg{3}{2}{-\frac{n-1}{2},-\frac{n}{2},\frac{a+n+1}{2}}
{\frac{3}{2},\;\frac{a-n+1}{2}-b}. 
\end{eqnarray*}
Here for each positive integer $n$ the two $\Gamma$-terms are
$\CC(a,b)$-multiples of the $\Gamma$-terms in (\ref{genkum}), so the
coefficients $P(n)$, $Q(n)$ are equal to $\CC(a,b)$-multiples of the
$\hpgo{3}{2}(1)$ series in this formula. But the correspondence depends
on whether $n$ is even or odd.}
formula (\ref{whip841}) to the right-hand side of 
the identity \cite[2.9(2)]{bateman}:
\begin{equation} \label{hpg21tr}
\hpg{2}{1}{a+n,\,b}{a-b}{-1} = 2^{-2b-n}\;\hpg{2}{1}{a-2b,\,-b-n}{a-b}{-1}.
\end{equation}
After this we sum up the terms with even and odd indexes separately,
transform the $\Gamma$-factors slightly and get formula (\ref{genkum})
with $P(n)$, $Q(n)$ defined by (\ref{whippn}-\ref{whipqn}).

Theorem \ref{gtheorem2} follows from Theorem \ref{gtheorem1}
by the following transformation of terminating $\hpgo{3}{2}(1)$
series (see \cite{specfaar}, proof of Cor.~3.3.4):
\begin{equation} \label{termhpg32}
\hpgg{3}{2}{-m,\,A,\,B}{E,\;F} =
\frac{(E\!-\!A)_m}{(E)_m}\,
\hpgg{3}{2}{-m,\;A,\;F\!-\!B}{1\!+\!A\!-\!E\!-\!m,\,F},
\end{equation}
where $m$ must be a non-negative integer. To make sure that
the interpretation of ill-defined hypergeometric series in
(\ref{pqdef1}) is correct for this transformation,
one may specialize $A$ to $-\nu/2$ or $-(\nu\!\pm\!1)/2$ with complex $\nu$
(instead of $-n/2$, etc.) and take the limit $\nu\!\to\!n$.

To prove Theorem \ref{gtheoremneg} we use Euler's integral again.
After rearranging the integrand in (\ref{rearrange}) as
$t^{A-1}(1\!-\!t)^{C-A+B-1}(1\!-\!t^2)^{-B}$ and expanding
$(1\!-\!t)^{C-A+B-1}$ we eventually get formula:
\begin{equation} \label{whip841a}
\hpg{2}{1}{\!A, B}{C}{-1}=
\frac{1}{2}\frac{\Gamma(C)\,\Gamma(1\!-\!B)}{\Gamma(A)\,\Gamma(C\!-\!A)}
\sum_{k=0}^{\infty} \frac{(\!A\!-\!B\!-\!C\!+\!1)_k}{k!}
\frac{\Gamma(\frac{A}{2}\!+\!\frac{k}{2})}
{\Gamma(\frac{A}{2}\!+\!\frac{k}{2}\!+\!1\!-\!B)}.
\end{equation}
Like in the proof of Theorem \ref{gtheorem1}, we apply this formula to
$\hpg{2}{1}{a-n-1,b}{a-b}{-1}$ transformed by (\ref{hpg21tr}), and add the
terms with even and odd indexes separately. The result is:
\begin{eqnarray} \label{whippn1}
P(-n-1) & = & 2^n\;\frac{(1-\frac{a}{2})_n}{(1-b)_n}\;
\hpgg{3}{2}{-\frac{n}{2},\,-\frac{n-1}{2},\,\frac{a}{2}-b}
{\frac{1}{2},\;\frac{a}{2}-n}, \\
\label{whipqn1} Q(-n-1) & = & -n\,2^n\;\frac{(\frac{1-a}{2})_n}{(1-b)_n}\;
\hpgg{3}{2}{-\frac{n-1}{2},\,-\frac{n-2}{2},\,\frac{a+1}{2}-b}
{\frac{3}{2},\;\frac{a+1}{2}-n}.
\end{eqnarray}
Comparing these expressions with (\ref{whippn}-\ref{whipqn}) gives Theorem
\ref{gtheoremneg}.
Q.E.D.\\

To get more expressions for $P(n)$ and $Q(n)$ one can use standard
transformations of terminating $\hpgo{3}{2}(1)$ series. For example, 
one may repeatedly apply 
(\ref{termhpg32}) or rewrite a terminating series in the reverse
order. In general, a terminating $\hpgo{3}{2}(1)$ series can be
transformed to 17 other terminating $\hpgo{3}{2}(1)$ series,
see \cite[sect.~8]{whipple24}, \cite[3.9]{bailey}. To give these
transformations a group structure one has to consider transpositions
of the two lower and two upper parameters as non-trivial transformations.
Then one gets a group of 72 elements which acts on the set of 18
hypergeometric series, see \cite{raojeugt}. The action of this group
can be summarized as follows. Let $y_0,\ldots,y_5$ be six
parameters satisfying $y_0+y_1+y_2=y_3+y_4+y_5=1-m$. Then the expression
\begin{equation} \label{revhpg32}
(y_0\!+\!y_4)_m\,(y_0\!+\!y_5)_m\;
\hpgg{3}{2}{-m,\,y_0\!+\!y_1\!-\!y_3,\,y_0\!+\!y_2\!-\!y_3}
{y_0\!+\!y_4,\;y_0\!+\!y_5}
\end{equation}
is invariant under the permutations within the sets $\{y_0,\,y_1,\,y_2\}$
and $\{y_3,\,y_4,\,y_5\}$, and gets multiplied by $(-1)^m$ when these two
sets are interchanged. For instance, formula (\ref{termhpg32}) corresponds
to the permutation $y_0\!\leftrightarrow\!y_5$, $y_1\!\leftrightarrow\!y_4$,
$y_2\!\leftrightarrow\!y_3$. 

Application of these transformations to the series (\ref{whippn}-\ref{whipqn})
or (\ref{pqdef1}) gives eight sets of 18 terminating $\hpgo{3}{2}(1)$ series,
one set for a choice of $P(n)$ or $Q(n)$, positive or negative and even or odd
$n$. The number of different hypergeometric series turns out to be 96. Here we
summarize a few interesting expressions for $n\ge 0$:
\begin{eqnarray} \label{altpn}
P(n) & = & \frac{\lfloor\frac{n}{2}\rfloor!}{2\cdot n!}\,
\left( {\textstyle \frac{1-a}{2}\!+\!b} \right)_{\lceil n/2 \rceil}\,
\hpgg{3}{2}{\!-\!\lceil\frac{n}{2}\rceil,
\frac{a+1}{2}\!+\!\lfloor\frac{n}{2}\rfloor,\frac{a}{2}\!-\!b}
{\frac{a}{2},\,\frac{a+1}{2}\!-\!\lceil\frac{n}{2}\rceil\!-\!b} \\
\label{altpn1} & = & \frac{1}{2^{n+1}}\,\frac{(b)_{\lceil n/2 \rceil}}
{(\frac{a}{2})_{\lceil n/2 \rceil}}\,\hpgg{3}{2}
{\!-\!\lceil\frac{n}{2}\rceil,1\!+\!\lfloor\frac{n}{2}\rfloor,
\frac{a}{2}-b}{\frac{1}{2},\;1\!-\!b\!-\!\lceil\frac{n}{2}\rceil},\\
\label{altqn} Q(n) & = & \frac{\lceil\frac{n}{2}\rceil!}{2\cdot n!}\,
\left( {\textstyle 1\!-\!\frac{a}{2}\!+\!b} \right)_{\lfloor n/2 \rfloor}\,
\hpgg{3}{2}{-\!\lfloor\frac{n}{2}\rfloor,
\frac{a}{2}\!+\!\lceil\frac{n}{2}\rceil,\frac{a+1}{2}\!-\!b}
{\frac{a+1}{2},\,\frac{a}{2}\!-\!\lfloor\frac{n}{2}\rfloor\!-\!b}\\
\label{altqn1} & = & \frac{n\!+\!1}{2^{n+1}}\,\frac{(b)_{\lfloor n/2 \rfloor}}
{(\frac{a+1}{2})_{\lfloor n/2 \rfloor}}\,\hpgg{3}{2}
{\!-\!\lfloor\frac{n}{2}\rfloor,1\!+\!\lceil\frac{n}{2}\rceil,
\frac{a+1}{2}-b}{\frac{3}{2},\;1\!-\!b\!-\!\lfloor\frac{n}{2}\rfloor}.
\end{eqnarray}
To get expressions for negative $n$ one may use Theorem
\ref{gtheoremneg}. Notice that series in (\ref{altpn1}) and (\ref{altqn1})
terminate for all $n$.
\\

In the rest of this section we follow Whipple's approach in \cite{whipple30},
where transformations of not necessarily terminating $\hpgo{3}{2}(1)$ series
are used to derive various identities with general $\hpgo{2}{1}(-1)$ series.
We concentrate on the $\hpgo{2}{1}(-1)$ series which are contiguous to the
series in Kummer's formula (\ref{ikummer}). Notice that proofs of Theorems
\ref{gtheorem1} and \ref{gtheoremneg} are valid for any complex values of $n$,
so that formula (\ref{genkum}) with $P(n)$ and $Q(n)$ defined by
(\ref{whippn}-\ref{whipqn}) or (\ref{whippn1}-\ref{whipqn1}) is true
for any complex $n$. Formula (\ref{genkum}) with $P(n)$, $Q(n)$ defined
by (\ref{pqdef1}) is also true for all $n$, see Whipple's formulas
(\ref{whip1}-\ref{whip2}) below. But one may check that in general
these $P(n)$ and $Q(n)$ are not the same.

Transformations of general $\hpgo{3}{2}(1)$ series were first derived by
Thomae, see \cite{thomae}. Whipple introduced notation
(see \cite{whipple24},\cite[3.5-7]{bailey}) which gives a group-theoretical
insight into those formulas. To begin with, there is an action of the
symmetric group $S_5$ on $\hpgo{3}{2}(1)$'s. Hardy described it in the notes
to lecture VII in \cite{hardyram} by saying that the function
\begin{equation} \label{invexpr}
\frac{1}{\Gamma(E)\,\Gamma(F)\,\Gamma(E\!+\!F\!-\!A\!-\!B\!-\!C)}
\;\hpgg{3}{2}{A,\,B,\,C}{E,\;F}
\end{equation}
is invariant under the permutations of $E$, $F$, $E\!+\!F\!-\!B\!-\!C$,
$E\!+\!F\!-\!A\!-\!C$ and $E\!+\!F\!-\!A\!-\!B$. For example, we have
(see \cite{specfaar}, Cor.~3.3.5):
\begin{equation} \label{hpg32tr}
\hpgg{3}{2}{A,\,B,\,C}{E,\;F} =
\frac{\Gamma(F)\,\Gamma(E\!+\!F\!-\!A\!-\!B\!-\!C)}
{\Gamma(F\!-\!C)\,\Gamma(E\!+\!F\!-\!A\!-\!B)}\,
\hpgg{3}{2}{E\!-\!A,\,E\!-\!B,\,C}{E,\,E\!+\!F\!-\!A\!-\!B}.
\end{equation}
An orbit of general $\hpgo{3}{2}(1)$ consists of 10 different series.
Note that the series in (\ref{invexpr}) converge when
$\re{E\!+\!F\!-\!A\!-\!B\!-\!C}>0$, and the whole expression
is well-defined and analytic for any parameters under this condition.
Function (\ref{invexpr}) can be analytically continued to the
region in the parameter space where at least one of the 10 series converges.

Further, a general $S_5$ orbit of $\hpgo{3}{2}(1)$'s is associated to
11 other orbits, so that we get sets of 120 {\em allied} $\hpgo{3}{2}(1)$
series, see \cite{whipple24}. For example\footnote{Whipple introduced for
$\hpgo{3}{2}(1)$ series six parameters $r_0,\ldots,r_5$ related by condition
\mbox{$\sum r_i\!=\!0$} so that: all allied series can be obtained by the
permutations of the six parameters and/or changing the sign of them all;
an $S_5$-orbit is determined by fixing a parameter and an element of the
set $\{+,-\}$, and $S_5$ permutes the remaining five parameters. Specifically,
one may choose that the $S_5$ action on (\ref{invexpr}) fixes $r_0$, and
take $E\!=\!1\!+\!r_4\!-\!r_0$.}, the series in (\ref{invexpr}) is allied to
\begin{equation}
\hpgg{3}{2}{A,\,1\!+\!A\!-\!E,\,1\!+\!A\!-\!F}{1\!+\!A\!-\!B,\,1\!+\!A\!-\!C}
\quad \mbox{and} \quad
\hpgg{3}{2}{E\!-\!A,\,E\!-\!B,\,E\!-\!C}{E,\,1\!+\!F\!-\!E}.
\end{equation}
In general, two allied series are not related by a two-term identity
like (\ref{hpg32tr}). But for any three allied series there is a linear
relation between them, with coefficients being $\Gamma$-terms. This also
gives three-term relations for the 12 functions of type (\ref{invexpr}),
and even defines their analytic continuation to the whole space of
parameters. Indeed, if the series in (\ref{invexpr}) diverges then its ally 
$\hpgg{3}{2}{1-A,\,1-B,\,1-C}{2-D,\;2-E}$ converges; for the third term
one can take convergent series from a similar pair of functions from other
$S_5$-orbits. Besides, all allied series converge in a neighborhood of
$A\!=\!B\!=\!C\!=\!1/2$, $E\!=\!F\!=\!1$.

In \cite{whipple30} Whipple applies the relations of allied series to a
general $\hpgo{2}{1}(-1)$ series by expressing it as $\hpgo{3}{2}(1)$ series
and considering it as a member of the corresponding allied family.
In particular, his formulas (3.1) and (3.51) read as:
\begin{eqnarray} \label{whip1}
\hpg{2}{1}{a+\nu,\,b}{a-b}{-1} & \!=\! &
\frac{\Gamma(a-b)\,\Gamma(\frac{a}{2})}{\Gamma(a)\,\Gamma(\frac{a}{2}-b)}\,
\hpgg{3}{2}{-\frac{\nu-1}{2},\,-\frac{\nu}{2},\,b}{-\nu,\,\frac{a+1}{2}}
\\ \label{whip2} & \!=\! &
\frac{\Gamma(a-b)\,\Gamma(\frac{a+1}{2})}{\Gamma(a)\,\Gamma(\frac{a+1}{2}-b)}
\,\hpgg{3}{2}{-\frac{\nu}{2},\,-\frac{\nu+1}{2},\,b}{-\nu,\,\frac{a}{2}}.
\end{eqnarray}
If $\nu\not\in\{0,1,2,\ldots,\}$ we may relate the $\hpgo{2}{1}(-1)$ series
with the $S_5$-orbit of the $\hpgo{3}{2}(1)$ series in
(\ref{whip1}-\ref{whip2}) and get many two- and three-term relations with
$\hpgo{2}{1}(-1)$ and $\hpgo{3}{2}(1)$ series. Some of these identities
make sense and are correct even if $\nu$ is a non-negative integer, because
singular $\Gamma$-factors cancel. For instance, formula (\ref{genkum})
with $P(n)$ and $Q(n)$ defined by (\ref{whippn}-\ref{whipqn}) is a three
term identity between allied series, see the last paragraph of
\cite{whipple30}.  Similarly, (potentially) terminating series in Whipple's
formulas (8.3) and (8.41) are derived from three term identities of allied
series.

On the other hand, the $\hpgo{3}{2}(1)$ series in (\ref{whip1}-\ref{whip2})
cannot be identified with the terminating series in the expressions in
(\ref{pqdef1}). One has to compute:
\begin{eqnarray*}
\lim_{\nu\to n}
\hpgg{3}{2}{\!-\frac{\nu}{2},\,-\frac{\nu+1}{2}, b}{-\nu,\,\frac{a}{2}}
= 2\,P(n)-\frac{1}{4^{n+1}}\,\frac{\poh{b}{n+1}}{\poh{\frac{a}{2}}{n+1}}
\,\hpgg{3}{2}{\frac{n+2}{2},\,\frac{n+1}{2},\,b\!+\!n\!+\!1}
{n\!+\!2,\,\frac{a}{2}\!+\!n\!+\!1},\\
\lim_{\nu\to n}
\hpgg{3}{2}{\!-\frac{\nu-1}{2}, -\frac{\nu}{2}, b}{-\nu,\,\frac{a+1}{2}} =
2\,Q(n)+\frac{1}{4^{n+1}} \frac{\poh{b}{n+1}}{\poh{\frac{a+1}{2}}{n+1}}
\hpgg{3}{2}{\frac{n+3}{2},\frac{n+2}{2},b\!+\!n\!+\!1}
{n\!+\!2,\,\frac{a+1}{2}\!+\!n\!+\!1}.
\end{eqnarray*}
In the sum of these two equalities the non-terminating $\hpgo{3}{2}(1)$ series
on the right-hand side cancel, since they are connected by transformation
(\ref{hpg32tr}). In this way identities (\ref{whip1}-\ref{whip2}) prove
Theorem \ref{gtheorem2}.

Moreover, the $\hpgo{3}{2}(1)$ series (\ref{whip1}-\ref{whip2}) can be
transformed by $S_5$ to four series which are well-defined and terminate
when $\nu$ is an (odd or even) positive integer $n$. Those terminating
series are presented in formulas (\ref{altpn}) and (\ref{altqn}). However,
this does not give expressions for $\hpg{2}{1}{\!a+n,b}{a-b}{-1}$ in terms
of one terminating $\hpgo{3}{2}(1)$ series, because the mentioned four
series diverge for $\nu>1/2$ (except when they terminate), and we cannot
use the $S_5$-invariance of the corresponding function in (\ref{invexpr}).
Notice, for example, that (\ref{hpg32tr}) implies a wrong relation
between the $\hpgo{3}{2}(1)$ series in (\ref{altpn}) and (\ref{altqn}).
As we see, Whipple's approach in \cite{whipple30} gets complicated in the
case $\nu$ in (\ref{whip1}-\ref{whip2}) is an integer,
and does not directly explain various expressions for our $P(n)$ and $Q(n)$. 

\section{A proof by Zeilberger's method}
\label{ctproof}

Here we prove Theorem \ref{gtheorem2} only. Theorem
\ref{gtheorem1} can be proved in the same way.

Let us define 
$S(n)=\hpg{2}{1}{\!a+n,\,b}{a-b}{-1}$.  The contiguity relation
\cite[2.8(28)]{bateman} between $\hpg{2}{1}{\!A+1,\,B}{C}{z}$,
$\hpg{2}{1}{\!A-1,\,B}{C}{z}$ and $\hpg{2}{1}{\!A,\,B}{C}{z}$
gives the following recurrence relation:
\begin{equation} \label{recrel}
2\,(n\!+\!a)\,S(n\!+\!1)-(3\,n\!+\!2\,a)\,S(n)+(n\!+\!b)\,S(n\!-\!1)=0.
\end{equation}
We claim that the sequences $P(n)$ and $Q(n)$ satisfy the same recurrence
relation. Following the ``creative telescoping'' method
of Zeilberger (\cite{zeilb, koepf}), let
\begin{equation}
p(n,k)=\frac{(-1)^k}{2\cdot 4^k}\,\frac{(n\!+\!1)\,(n\!-\!k)!}
{(n\!-\!2k\!+\!1)!\,k!}\,\frac{(b)_k}{(\frac{a}{2})_k}
\end{equation}
be the $k$th summand of $P(n)$ in (\ref{pqdef1}). We set $p(n,k)=0$ for
$k>\lceil n/2 \rceil$. Also define
\[
r_1(n,k)=-\frac{2\,k\,(n\!-\!k\!+\!1)\,(a\!+\!2k\!-\!2)}
{(n\!-\!2k\!+\!2)\,(n\!-\!2k\!+\!3)},\qquad
R_1(n,k)=r_1(n,k)\,p(n,k).
\]
One can check that
\begin{eqnarray*}
2\,(n\!+\!a)\,p(n\!+\!1,k)\!-\!(3n\!+\!2a)\,p(n,k)\!+\!
(n\!+\!b)\,p(n\!-\!1,k) \!=\! R_1(n,k\!+\!1)\!-\!R_1(n,k),
\end{eqnarray*}
so
\begin{eqnarray} \label{rdef1}
2(n\!+\!a)\,P(n\!+\!1)-(3\,n\!+\!2\,a)\,P(n)+(n\!+\!b)\,P(n\!-\!1)
 & \!\!=\!\! \nonumber \\
\sum_{k=0}^{\lfloor n/2 \rfloor} \big( R_1(n,k\!+\!1)-R_1(n,k) \big)
-R_1(n,\lceil\textstyle\frac{n+1}{2}\rceil) & \!\!=\!\! & 0.
\end{eqnarray}
Although this looks like an artificial trick, we follow the standard
Wilf-Zeilberger method of proving combinatorial identities, see
\cite{zeilb, koepf}. The expression $r_1(n,k)$ is the {\em certificate}
of our standardized proof. Given $p(n,k)$ the recurrence relation for
$P(n)$ and the certificate $r_1(n,k)$ can be found by Zeilberger's
algorithm. 
This algorithm is implemented in computer algebra
packages {\sf EKHAD} (see \cite{ekhad}, command {\sf ct}) and {\sf hsum.mpl}
(see \cite{pkoepf}, command {\sf sumrecursion} with option
{\sf certificate=true}). Also check \cite{koornvid} for a link to a
{\sf Maple} worksheet for this proof.
The finite sums in this proof require some attention, since they are not
natural according to \cite{koepf}.

In the same way: 
\begin{eqnarray} \label{rdef2}
2\,(n\!+\!a)\,Q(n\!+\!1)-(3\,n\!+\!2\,a)\,Q(n)+(n\!+\!b)\,Q(n-1) & 
\!\!=\!\! \nonumber \\
\sum_{k=0}^{\lfloor (n-1)/2 \rfloor} \big(
R_2(n,k\!+\!1)-R_2(n,k) \big)-R_2(n,\lceil\textstyle\frac{n}{2}\rceil)
& \!\!=\!\! & 0,
\end{eqnarray}
where
\begin{equation}
R_2(n,k) = \frac{2\,k\,(n\!-\!k\!+\!1)\,(a\!+\!2k\!-\!1)}
{(n\!-\!2k\!+\!1)\,(n\!-\!2k\!+\!2)} \cdot
\frac{(-1)^k}{2\cdot4^k}\,{n-k \choose k}\,\frac{(b)_k}{(\frac{a+1}{2})_k}.
\end{equation}
is the $k$th summand of $Q(n)$ in (\ref{pqdef1}) multiplied by the
corresponding certificate.

Note that the condition $(a)_n\neq 0$ ensures that recurrence relation
(\ref{recrel}) does not degenerate to a first order relation until we evaluate
$P(n)$ and $Q(n)$. It remains to check that formula (\ref{genkum}) holds for
two initial values of $n$. Kummer's identity (\ref{ikummer}) suggests
$P(-1)\!=\!1$ and $Q(-1)\!=\!0$, which fits into the recurrence relation.
We may use Gauss' contiguity relation \cite[2.8(38)]{bateman} between
$\hpg{2}{1}{\!A,\,B}{C+1}{z}$, $\hpg{2}{1}{\!A,\,B}{C}{z}$ and
$\hpg{2}{1}{\!A-1,\,B}{C}{z}$ to obtain
\begin{equation}
(a-2b)\,\frac{\Gamma(1+a-b)\,\Gamma(1+\frac{a}{2})}
{\Gamma(1+a)\,\Gamma(1+\frac{a}{2}-b)}-2\,(a-b)\,S(0)+(a-b)\,S(-1)=0.
\end{equation}
This implies the correct $P(0)\!=\!1/2$ and $Q(0)\!=\!1/2$ and completes the
proof.

Note that the Gauss contiguity relations hold for analytic extensions of
hypergeometric functions on $\CC\setminus [1,\infty)$. Therefore this proof
does not require convergence of the $\hpgo{2}{1}(-1)$ series.
Q.E.D.\\

In fact, sequences $P(n)$ and $Q(n)$ satisfy recurrence relation (\ref{recrel})
for all $n$. The recurrence can be directly verified for $n=-2,-1,0$.
The values of $P(n)$ and $Q(n)$ for $n=-3,-2,-1,0,1$ are
\[ \begin{array}{cccccl}
\frac{2\,(a-2)\,(a-b-2)}{(b-1)\,(b-2)}, & \frac{a-2}{b-1}, &
1, & \frac{1}{2}, & \frac{a-b}{2\,a} & \mbox{and} \vspace{2pt} \\
-\frac{2\,(a-1)\,(a-3)}{(b-1)\,(b-2)}, & -\frac{a-1}{b-1}, &
0, & \frac{1}{2}, & \frac{1}{2} & \mbox{respectively.}
\end{array}
\]
To compute the same recurrence relation for negative $n$ one can use Theorem 1.
Alternatively, one may choose an expression for $P(n)$ and $Q(n)$ for negative
$n$, say (\ref{whippn1}-\ref{whipqn1}), and compute the recurrence relation
with Zeilberger's algorithm.

To show equalities like in (\ref{altpn}-\ref{altqn}) by Zeilberger's method
one would have to compute the recurrences for odd and even integers separately.
Recurrence relation (\ref{recrel}) for any such expression and for all $n$ can
be computed using contiguity relations for $\hpgo{3}{2}(1)$ series.
As it is known (see \cite[3.7]{specfaar}), contiguous $\hpgo{3}{2}(1)$ series
satisfy three-term relations (with coefficients being rational functions in
the parameters of those series), just like contiguous $\hpgo{2}{1}(z)$ series.

\section{Algorithmic aspects}
\label{comalge}

The generalized formula (\ref{genkum}) can be specialized so that $P(n)$
or $Q(n)$ vanishes, giving an evaluation of $\hpgo{2}{1}(-1)$ series with
a single $\Gamma$-term. For example, 
\begin{equation}
Q(-4)=-4\,\frac{(a-1)\,(a-3)\,(2a-b-7)}{(b-1)\,(b-2)\,(b-3)},
\end{equation}
so if $b=2a-7$ then $Q(-4)=0$, which implies
\begin{equation} \label{specfo1}
\hpg{2}{1}{3-c,\,7-2\,c}{c}{-1}=\frac{3}{4}\,
\frac{\Gamma(c)\,\Gamma(3-\frac{c}{2})}
{\Gamma(5-c)\,\Gamma(\frac{3\,c}{2}-2)}.
\end{equation}
Further, $P(-5)=0$ if $2\,a^2-4\,a\,b+b^2-12\,a+17\,b+12=0$. Parameterizing
the curve given by this equation we get
\begin{equation} \label{specfo2}
\hpg{2}{1}{-\frac{2t^2-7t+6}{t^2-2},\,\frac{t^2+4t-8}{t^2-2}}
{\frac{2t^2+3t-8}{t^2-2}}{-1}=\frac{t^2+3\,t-6}{t\,(t-1)}\,
\frac{\Gamma\left(\frac{3t-4}{t^2-2}\right)\,
\Gamma\left(\frac{t^2+7t-12}{2\,(t^2-2)}\right)}
{\Gamma\left(\frac{7t-10}{t^2-2}\right)\,
\Gamma\left(\frac{t\,(t-1)}{2\,(t^2-2)}\right)}.
\end{equation}

It could be expected that formulas like (\ref{specfo1}) can be proved
automatically by current computer algebra algorithms, say
by Wilf-Zeilberger method. As it is demonstrated in \cite{koornid},
this method or Zeilberger's algorithm can be adapted to non-terminating
hypergeometric series if one can justify the ``creative telescoping''
trick by dominated convergence, and the hypergeometric series can be
evaluated in the limit $n\to\infty$, where $n$ is a discrete parameter.
In general non-terminating hypergeometric series is given without a
discrete parameter, so it must be introduced by an algorithm.
For example, after substitution $a\mapsto a+2n$ one can prove Kummer's
formula (\ref{ikummer}) with Wilf-Zeilberger method, see
\cite{gauthier}.

In the case of equation (\ref{specfo1}) we may substitute
$c\mapsto c+n$ and apply Zeilberger's algorithm to get the right first order
difference equation. However, we cannot evaluate the hypergeometric series
neither in the limit $n\to\infty$, nor for a finite value of $n$.
What we can do is to 
combine explicitly Gauss' contiguity relations in such a way that we
``accidentally'' get a two-term relation where one of the terms can
be evaluated by Kummer's formula. For example, the relation between
contiguous $\hpg{2}{1}{A,\,B}{C}{z}$, $\hpg{2}{1}{A+1,\,B-2}{C}{z}$ and,
say, $\hpg{2}{1}{A,\,B-1}{C}{z}$, after the specialization
$(A,\,B,\,C,\,z)\mapsto (3\!-\!c,\,7\!-\!2c,\,c\!-\!1)$ becomes
\begin{equation}
\hpg{2}{1}{3-c,\,7-2\,c}{c}{-1}=
\frac{3}{4}\,\hpg{2}{1}{4-c,\,5-2\,c}{c}{-1},
\end{equation}
In this way even the exotic (\ref{specfo2}) can be proved.

This shows that relations between contiguous hypergeometric series can be
useful for finding new ``non-standard'' evaluations of $\hpgo{2}{1}$ series.
One may take such a relation and try to find families of its two term
specializations with a discrete parameter $n$. This would give a first order
recurrence relation, and if the series can be evaluated in the limit
$n\to\infty$ one gets a (perhaps) new formula! Relations between contiguous
series also give a way to compute recurrence relation, alternative to
Zeilberger's algorithm.

In \cite{koornvid} there is a link to {\sf Maple} routines which for given
three integer vectors $(k_i,l_i,m_i)$ for $i=1,2,3$ derive a
$\CC(A,B,C,z)$-linear relation between three contiguous functions
$\hpg{2}{1}{A+k_i,\,B+l_i}{C+m_i}{z}$. Computer experiments found many
first order recurrence relations for some values $z=1/4$, $1/3$, $1/9$,
$\exp(i\pi/3)$, $3-2\sqrt{2},\ldots$, some of them can be successfully
solved. It is an interesting question which $\hpgo{2}{1}(z)$ series can be
evaluated in terms of $\Gamma$-function. So far produced evaluations can be
obtained using quadratic or cubic transformations.
\\

Here we generalize a few known formulas of the same
type as (\ref{genkum}). They were obtained by considering relations between
three contiguous hypergeometric series where two of them can be evaluated
by a known formula, and trying to express the coefficients in these relations
as hypergeometric series. This was done by considering partial fraction
decomposition of those coefficients empirically. The formulas can be proved
by showing that all three terms in a formula satisfy the same recurrence
relation by Zeilberger's algorithm, and checking the identity for a couple
of values of the discrete parameter.

We start with a generalization of Gosper's ``non-standard'' evaluations of
$\hpgo{2}{1}(1/4)$ series, see \cite[1/4.1-2]{gosper4th}. A generalization is
\begin{eqnarray} \label{gengosper}
\hpg{2}{1}{-a,\,\frac{1}{2}}{2a\!+\!\frac{3}{2}\!+\!n}{\frac{1}{4}} =
\frac{2^{n+3/2}}{3^{n+1}}
\frac{\Gamma(a\!+\!\frac{5}{4}\!+\!\frac{n}{2})\,
\Gamma(a\!+\!\frac{3}{4}\!+\!\frac{n}{2})\,\Gamma(a\!+\!\frac{1}{2})}
{\Gamma(a\!+\!\frac{7}{6}\!+\!\frac{n}{3})\,
\Gamma(a\!+\!\frac{5}{6}\!+\!\frac{n}{3})\,
\Gamma(a\!+\!\frac{1}{2}\!+\!\frac{n}{3})}\,K(n) \nonumber \\
-(-3)^{n-2}\,2^{3/2}\,\frac{\Gamma(a\!+\!\frac{5}{4}\!+\!\frac{n}{2})\,
\Gamma(a\!+\!\frac{3}{4}\!+\!\frac{n}{2})\,\Gamma(a\!+\!1)}
{\Gamma(a\!+\!\frac{3}{2})\,\Gamma(a\!+\!\frac{1}{2}\!+\!\frac{n}{2})\,
\Gamma(a\!+\!1\!+\!\frac{n}{2})}\,L(n),
\end{eqnarray}
where 
\[ K(1)=L(0)=0, \qquad K(0)=L(1)=1, \]
for $n>1$:
\begin{eqnarray*}
K(n) & \!\!=\!\! & (-1)^n\,\sum_{k=\lceil n/3\rceil}^{\lfloor n/2\rfloor}
\frac{27^k}{4^k}\,\frac{n\,(k\!-\!1)!}{(n\!-\!2k)!\,(3k\!-\!n)!}\,
\frac{(a\!+\!\frac{1}{2})_k}{(a\!+\!1)_k},\\
L(n) & \!\!=\!\! & \hpgg{4}{3}{-\frac{n-1}{3},\,-\frac{n-2}{3},\,
-\frac{n-3}{3},\,a\!+\!1}{-\frac{n-2}{2},\,-\frac{n-3}{2},\,
a\!+\!\frac{3}{2}},  
\end{eqnarray*}
and for $-n<0$:
\begin{eqnarray*}
K(-n) & \!\!\!\!=\!\!\!\! &
\hpgg{4}{3}{\!-\frac{n}{3}, -\frac{n-1}{3}, -\frac{n-2}{3}, -a}
{-\frac{n-1}{2}, -\frac{n-2}{2}, -a\!+\!\frac{1}{2}}
= \sum_{k=0}^{\lfloor n/3\rfloor}
\frac{(-4)^k}{27^k} \frac{n\,(n\!-\!2k\!-\!1)!}{(n\!-\!3k)!\,k!}
\frac{(-a)_k}{(-a\!+\!\frac{1}{2})_k},\\
L(-n) & \!\!\!\!=\!\!\!\! &
(-1)^n \sum_{k=\lceil(n+1)/3\rceil}^{\lfloor(n+1)/2\rfloor}
\frac{27^k}{4^k}\,\frac{(n\!+\!1)\,(k\!-\!1)!}
{(n\!-\!2k\!+\!1)!\,(3k\!-\!n\!-\!1)!}\,
\frac{(-a\!-\!\frac{1}{2})_k}{(-a)_k}.
\end{eqnarray*}
Gosper has found the special cases $n=0,1$. The $\Gamma$-factors to
$K(n)$ and $L(n)$ are $\CC(a)$-multiples of these two Gosper's evaluations
(respectively) for each $n$. All three terms in
(\ref{gengosper}) satisfy the recurrence relation
\begin{eqnarray*}
2\,(n\!+\!2a\!+\!1)\,(2n\!+\!6a\!+\!3)\,S(n\!+\!1)+
(2n\!+\!4a\!+\!3)\,(4n\!+\!6a\!+\!1)\,S(n) \\
-3\,(2n\!+\!4a\!+\!1)\,(2n\!+\!4a\!+\!3)\,S(n\!-\!1) & \!\!=\!\! & 0.
\end{eqnarray*}

Next we recall the classical Dixon's identity which evaluates well-poised
$\hpgo{3}{2}(1)$ series, see \cite[3.1]{bailey}. We generalize it as follows:
\begin{eqnarray} \label{gendixon}
\hpgg{3}{2}{a\!+\!n,\,b,\,c}{a\!-\!b,\,a\!-\!c}
& \!=\! & \frac{\tilde{P}(n)}{2}\,
\frac{\Gamma(\frac{a+1}{2})\,\Gamma(a\!-\!b)\,\Gamma(a\!-\!c)\,
\Gamma(\frac{a+1}{2}\!-\!b\!-\!c)}{\Gamma(a)\,\Gamma(\frac{a+1}{2}\!-\!b)
\,\Gamma(\frac{a+1}{2}\!-\!c)\,\Gamma(a\!-\!b\!-\!c)} \nonumber \\
& & + \frac{\tilde{Q}(n)}{2}\,\frac{\Gamma(\frac{a}{2})\,\Gamma(a\!-\!b)\,
\Gamma(a\!-\!c)\,\Gamma(\frac{a}{2}\!-\!b\!-\!c)}{\Gamma(a)\,
\Gamma(\frac{a}{2}\!-\!b)\,\Gamma(\frac{a}{2}\!-\!c)\,\Gamma(a\!-\!b\!-\!c)},
\end{eqnarray}
where $\tilde{P}(-1)=1$, $\tilde{Q}(-1)=0$, then
for $n\ge 0$:
\begin{eqnarray*}
\tilde{P}(n) = \hpgg{4}{3}{-\frac{n}{2},\,-\frac{n+1}{2},\,b,\,c}
{-n,\,\frac{a}{2},\,\frac{1-a}{2}\!+\!b\!+\!c},\quad
\tilde{Q}(n) = \hpgg{4}{3}{-\frac{n-1}{2},\,-\frac{n}{2},\,b,\,c}
{-n,\,\frac{1+a}{2},\,1\!-\!\frac{a}{2}\!+\!b\!+\!c},
\end{eqnarray*}
and for $-n<0$:
\begin{eqnarray*}
\tilde{P}(-n\!-\!1) \!\!& = &\!\!
2^{2n}\,\frac{(1\!-\!\frac{a}{2})_n\,(\frac{1+a}{2}\!-\!b\!-\!c)_n}
{(1\!-\!b)_n\,(1\!-\!c)_n}\,
\hpgg{4}{3}{-\frac{n}{2},\,-\frac{n-1}{2},\,b\!-\!n,\,c\!-\!n}
{1\!-\!n,\,\frac{a}{2}\!-\!n,\,\frac{1-a}{2}\!+\!b\!+\!c\!-\!n},\\
\label{pqalt4} \tilde{Q}(-n\!-\!1) \!\!& = &\!\!
-2^{2n}\,\frac{(\frac{1-a}{2})_n\,(\frac{a}{2}\!-\!b\!-\!c)_n}
{(1\!-\!b)_n\,(1\!-\!c)_n}\,
\hpgg{4}{3}{-\frac{n-1}{2},\,-\frac{n-2}{2},\,b\!-\!n,\,c\!-\!n}
{1\!-\!n,\,\frac{1+a}{2}\!-\!n,\,1\!-\!\frac{a}{2}\!+\!b\!+\!c\!-\!n}.
\end{eqnarray*}
Dixon's identity is the special case $n=-1$.  This generalized formula is
a relation between contiguous $\hpgo{3}{2}(1)$ series in explicit form.
For positive $n$ it is strikingly similar to generalization
(\ref{genkum},\ref{pqdef1}) of Kummer's identity. In fact, the generalization
in Theorem \ref{gtheorem2} is the limiting case $c\to\infty$ of
(\ref{gendixon}), just as Kummer's formula is the limiting case of
Dixon's identity. The recurrence relation for the three terms in
(\ref{gendixon}) is:
\begin{eqnarray*}
(n\!+\!a)\,(n\!-\!a\!+\!2b\!+\!2c\!+\!1)\,S(n\!+\!1)+
(n\!+\!b)\,(n\!+\!c)\,S(n\!-\!1) \nonumber \\
-(2n^2\!+\!3bn\!+\!3cn\!+\!n\!-\!a^2\!+\!2ab\!+\!2ac\!+\!a)\,S(n)
& \!=\! & 0.
\end{eqnarray*}

More evaluations of the same type can be obtained using standard
transformations of $\hpgo{2}{1}(z)$ series to $\hpgo{2}{1}(z/(z-1))$
series, see \cite[2.9(3-4)]{bateman}.
Applying them to the generalized Kummer's formula (\ref{genkum}) gives
evaluations of $\hpgo{2}{1}(1/2)$ which generalize classical formulas
of Gauss and Bailey, see \cite[2.4]{bailey}. The same transformation
of (\ref{gengosper}) gives evaluation of $\hpgo{2}{1}(-1/3)$. Similarly,
one can apply (\ref{hpg32tr}) to identity (\ref{gendixon}) and get
generalizations of Watson's and Whipple's formulas \cite[3.3-4]{bailey}.

All these formulas evaluate hypergeometric series which are contiguous
to a series which evaluation is known. In order to find these formulas
automatically one needs an algorithm which would find the solutions of
a recurrence relation in form of terminating hypergeometric series.

\bibliographystyle{alpha}
\bibliography{hypergeometric}

\end{document}